\documentclass[12pt, a4paper, reqno]{amsart}

\usepackage[utf8]{inputenc}
\usepackage{amsfonts,amssymb,amsxtra,amsthm,amsmath,amscd,mathrsfs,cite}
\usepackage[all]{xy}
\usepackage{enumerate}  
\usepackage{tikz}
\usepackage[normalem]{ulem}
\usepackage{soul}
\usepackage{color}

\newtheorem{theorem}{Theorem}[section]
\newtheorem{lemma}[theorem]{Lemma}

	\theoremstyle{remark}
	
	\numberwithin{equation}{section}
	\addtocounter{footnote}{1}

	\title[M\"obius disjointness conjecture   for Affine linear flow in short intervals]  
	{M\"obius disjointness conjecture \\ for Affine linear flow in short intervals} 
	
	\author{Qingyang Liu, Jing Ma}
	
	\address{School of Mathematics, Renmin University of China, Beijing, 100872, People's Republic of China}
	\email{qingyangliu@ruc.edu.cn}
	\address{School of Mathematics, Jilin University, Changchun, 130012, People's Republic of China}
	\email{jma@jlu.edu.cn}


\begin{document}
		\begin{abstract}
			In this paper we will generalize Liu and Sarnak's work on the M\"obius Disjointness Conjecture of 
			affine linear maps on a compact abelian group to short interval.
			The main result of this paper is that the M\"{o}bius Disjointness Conjecture of Sarnak 
			holds in short intervals $(N-M, N]$ with 
			$N^{5/8+\varepsilon}\leqslant  M\leqslant  N^{1-\varepsilon}$ for the affine linear flow  $(T, X)$ given by $T(x) = Ax + b$,
			where $A$ is an automorphism of $X$ and $b \in X$, and for the special case $X = \mathbb{T}^2$ we can also give explicit convergence rates.
		\end{abstract}
		
		\maketitle
		
		\section{Introduction}
		\subsection{The M\"{o}bius Disjointness Conjecture of Sarnak} 	
		
		Let $\mu(n)$ denote the M\"{o}bius function, defined by $\mu(1)=1$, $\mu(n)=(-1)^k$ if $n$ is a product of 
		$k$ distinct primes, and $\mu(n)=0$ otherwise. 
		The behavior of $\mu(n)$ is of fundamental importance in analytic number theory, 
		particularly in the study of prime numbers. 
		A flow is a pair $(X, T)$, where $X$ is a compact topological space and 
		$T \colon X \to X$ is a continuous map.  
		Sarnak's M\"{o}bius Disjointness Conjecture \cite{Sar12}, MDC in brief, states the following.
		
		\medskip
		
		\noindent
		\textbf{M\"{o}bius Disjointness Conjecture (Sarnak).}  
		\emph{The M\"{o}bius function $\mu$ is linearly disjoint from every zero-entropy flow $(X,T)$.  
			Equivalently, for any $f \in C(X)$ and any $x \in X$,
		\begin{equation}\label{MDC}
				\lim_{N\to\infty} \frac1N \sum_{n \leqslant  N} \mu(n) \, f(T^n x) = 0.
		\end{equation}}
		
		\medskip
		
		This conjecture establishes a profound connection between number theory and dynamical systems theory, 
		stimulating extensive research at the interface of analytic number theory and ergodic theory. 
		The simplest instance occurs when $X$ is a single point set and $T$ the identity map; this case is equivalent 
		to the Prime Number Theorem.  
		The first nontrivial case, where $X = \mathbb{T}$ is the unit circle and $T$ a translation, 
		corresponds to Vinogradov's estimate for exponential sums over primes, which played a crucial role in 
		his proof of the ternary Goldbach problem.

		Over the years, the conjecture has been established for a number of flows, 
		see for example \cite{Bou13, BouSarZie13, Sar12, LiuSar15, GreTao12,FerKulLem18} and the references therein. %, and most of these flows are regular, which means that the limit of their ergodic averages
		%\begin{eqnarray}\label{def/Birkh}
		%	\lim_{N\rightarrow\infty}\frac{1}{N}\sum_{n\leqslant  N} f(T^n x)
		%\end{eqnarray}
		%exist for every $f\in C(X)$ and every $x\in X$. 
	 
	%Other related work...... 
	 
Liu and Sarnak \cite{LiuSar15} investigated a class of regular dynamical systems defined by affine linear transformations acting on compact abelian groups. For a compact abelian group $X$, an affine linear map $T\colon X\to X$ takes the canonical form
	\begin{equation}\label{eqT}
		T(x) = Ax + b,
	\end{equation} 
where $A$ denotes a group automorphism of $X$ and $b\in X$ is a constant translation element; readers may refer to  \cite{Hahn1963affine,HoareParry1966quasi} for the foundational theory of affine flows. As a central theorem in their work, Liu and Sarnak verified the  M\"{o}bius Disjointness Conjecture for all zero-entropy affine linear flows \(\mathcal{X}=(T,X)\). 	
		
In this paper we will generalize Liu and Sarnak's work in \cite{LiuSar15} to short interval. %, and generalize the automorphism $A$ of $X$ to an homomorphism of $X$ such that 
%$$
%A(x,y)=Ax+Ay, \		\qquad \forall x,y\in X.
%$$ 		
		
\begin{theorem}\label{thm1}
	Let $\mathscr{X} = (X, T)$ be an afﬁne linear ﬂow on a compact abelian group $X$ which is of zero 	entropy, where  $T$ act on $X$ by the affine linear map \eqref{eqT},
	where $A$ is a group automorphism of $X$ and $b\in X$ is a constant translation element. 
	Then, $\mu$ is linearly disjoint from $\mathscr{X}$ in the short interval $N^{5/8+\varepsilon}\leqslant  M\leqslant  N^{1-\varepsilon}$ with  $N\geqslant 3$, i.e., 
	\begin{equation}\label{thm1goal}		
			\lim_{N\to\infty} \frac1M \sum_{N<n\leqslant  N+M} \mu(n) \, f(T^n x) = 0, \qquad\forall x\in X, \ f\in C(X).	
	\end{equation}
\end{theorem}		

For the special case $X = \mathbb{T}^2$, we can also derive a short-interval version of the  M\"{o}bius Disjointness Conjecture  with explicit convergence rates.

\begin{theorem}\label{thm2}
	Let $\mathscr{X} = (\mathbb{T}^2, T)$ be an  affine linear flow  of zero entropy with  $T$ defined by the affine linear map \eqref{eqT},
	where $A$ is a group automorphism of $\mathbb{T}^2$ and $b=(b_1,b_2)\in \mathbb{T}^2$. 
	Then, for $N^{5/8+\varepsilon}\leqslant  M\leqslant  N^{1-\varepsilon}$ with  $N\geqslant 3$ and $B>0$, we have
	\begin{equation}\label{thm1goal}		
		\frac1M \sum_{N<n\leqslant  N+M} \mu(n) \, f(T^n x) 
		\ll \log^{-B}N, \qquad\forall x\in X.	
	\end{equation}
	for any Lipschitz function $f$ on $\mathbb{T}^2$
\end{theorem}		

\medskip 

\noindent\textbf{Notation.}  
We use standard notation in number theory. For example ${\textup{e}}(x) = e^{2\pi i x}$.  The Landau-Vinogradov symbols $\ll, \gg$, and $O$ are used with their standard meanings. Throughout this paper, $\varepsilon$ always denotes an arbitrarily small positive constant.
The set of integers is denoted by $\mathbb{Z}$ and  the set of real number is denoted by $\mathbb{R}$. We write $\mathbb{T}$ for the quotient group $\mathbb{R}/\mathbb{Z}$.

\section{Preliminaries}

To prove Theorem \ref{thm1} and Theorem \ref{thm2}, we will need the following number-theoretic result, which is a direct corollary of \cite[Corollary 1.3]{TerenceTao2023Higheruniformity} by Matom{\"a}ki et al. Its proof follows along the same lines as the argument given in \cite{TerenceTao2023Higheruniformity} for the case   $\kappa=1$.

\begin{lemma}\label{youlideguji}
	\cite[Corollary 1.3]{TerenceTao2023Higheruniformity} 
	Let $\varepsilon>0$, $0\leqslant  l< q$,  $N^{5/8+\varepsilon}\leqslant  M\leqslant  N^{1-\varepsilon}$ with $N\geqslant 3$, 
and $P:\mathbb{Z}\rightarrow \mathbb{R}$ be a polynomial of degree $d\geqslant 1$. 
	Then, for any $B>0$ we have
	\[
	\sum_{\substack{N<n\leqslant  N+M \\n \equiv l \pmod{q}}}\mu(n){\textup{e}}(P(n))\ll_{d, B, \varepsilon, q} \frac{M}{\log^{B} N}.
	\]
\end{lemma}

In the proof of Theorem \ref{thm2}, we will also need the following lemma on topological entropys which  is a special case of \cite[Theorem 8.15]{Wal82}.

\begin{lemma}\label{lemh}	 
	Suppose $S: \mathbb{T}^m \to  \mathbb{T}^m$ is an affine transformation,  $S(x) = a\cdot A(x)$, where $a\in \mathbb{T}^m$ and $A$ is a surjective endomorphism of $\mathbb{T}^m$. Then
	\[
	h_{\textup{top}}(S)=h_{\textup{top}}(A)
	= \sum_{|\lambda_i|>1} \log|\lambda_i|,
	\]
	where $\lambda_1,\dots,\lambda_m$ are the eigenvalues of the matrix which represents $A$.
\end{lemma}

%The second lemma we need in the proof of Theorem \ref{thm2} is a  result of Zhan \cite[Theorem 5]{Zha91}  on exponential sum over primes in short intervals. In \cite[Theorem 5]{Zha91} only the case $q=r=1$ is treated, but it is plain to see that his proof carries over to the general case below. 
%\begin{lemma}\label{zhan}
%	Let $q$ and $r\leqslant   q$ be fixed co-prime positive integers. 	Let $B>0$ and $\varepsilon>0$ be arbitrarily large and small respectively, 	and assume $N^{5/8+\varepsilon}\leqslant  M\leqslant  N$. Then 
%	\[
%	\sum_{\substack{N-M<n\leqslant  N\\ n\equiv r \bmod q}} 
%	\mu(n){\textup{e}}(\alpha n)\ll_{B,\varepsilon, q} \frac{M}{\log^B N}  
%	\]
%	holds uniformly for $\alpha\in \mathbb{R}$.  
%\end{lemma}

%%%%%%%%%%%%%%%%%%%%%%%%%%%%%%%%%%%%%%%%%%%%%%%%%%% 
 
\section{Proof of Theorem \ref{thm1}}
\subsection{Reduction to the toral case.}

We first reduce the theorem to the case that $X$ is a torus, %(not necessarily connected), 
i.e., $X = \mathbb{T}^r \times C$ for some integer $r \geqslant 0$ and finite abelian group $C$. 
	
Since the linear combinations of characters in the (discrete) dual group $ \widehat{X}$ of $X$ are dense in $C(X)$, it suffices to show that for every $x \in X$ and $\psi \in  \widehat{X} $ we have
	\begin{equation}\label{thm1goal}
	\frac{1}{M} \sum_{N<n\leqslant  N+M} \mu(n) \psi(T^n x) \to 0, \quad \text{as } N \to \infty. 
	\end{equation} 
	
Define the act of $A$ on $\widehat{X}$ by 
	$$
	A\phi(x) := \phi(Ax), \qquad   \phi \in \widehat{X},\ x \in X.
	$$
Then $A\phi \in \widehat{X}$ for each $\phi \in \widehat{X}$ since $A$ is an automorphism of $X$, i.e., $\widehat{X}$ is $A$-invariant.
Noticing that the system $(X,T)$ is distal (from \cite{LiuSar15} or \cite{Lian17}), applying Aoki's result in \cite[page 13]{Aoki} we conclude  that $C_\psi$ is finitely generated, where for given $\psi \in \widehat{X} $, $C_\psi$ denotes the smallest  
$A$-invariant closed subgroup of $\widehat{X}$ containing $\psi$. For any subgroup $C$ of $ \widehat{X} $, define the annihilator of $C$ as 
	$$
	C^\perp = \{ x \in X \mid c(x) = 1 \text{ for all } c \in C \}.
	$$
Then, $C^\perp$ is a closed subgroup of $X$. In particular, for $C=C_\psi$, $X/C_\psi^\perp$ is a compact quotient group with $\widehat{X/C_\psi^\perp} = C_\psi$ finitely generated.  %This in turn implies that $\widehat{X/C_\psi^\perp}  = C_\psi$ is finitely generated. 
Therefore, as the dual of the finitely generated group $\widehat{X/C_\psi^\perp}$,  $X/C_\psi^\perp$ is isomorphic to $\mathbb{T}^r \times C$ for some $r \geqslant 0$ and finite abelian group $C$.  
	
Under these notations, for any  $y \in C_\psi^\perp$, we have $\phi(y)=1$ for any $\phi\in C_\psi$. Then, noting that $C_\psi$ is $A$-invariant, 	we have $A\phi\subseteq C_\psi$ which implies $\phi(Ay)=A\phi(y)=1$ for all $\phi\in C_\psi$ and $y\in  C_\psi^\perp$, which implies $Ay\in C_\psi^\perp$ for $y\in C_\psi^\perp$ (i.e., $C_\psi^\perp$ is closed under the action of $A$) and then
\[
T(x+y) = A(x+y)+b = Ax+b + Ay \equiv Ax+b \ (\bmod \ C_\psi^\perp), \quad \forall x \in X.
\]	
Hence $T(x+y) = T(x) \ (\bmod \ C_\psi^\perp)$, i.e., $T$ induces an affine linear map $T_\psi$ of $X/C_\psi^\perp$. Put another way, the flow $\mathcal{X}_\psi = (T_\psi, X/C_\psi^\perp)$ is a factor of $\mathcal{X} = (T, X)$ for any  $\psi \in \widehat{X} $. Since we are assuming that $\mathcal{X}$ has zero entropy it follows that so does $\mathcal{X}_\psi$. 

Moreover, for $n \geqslant 0$
\[
\psi(T^n x) = \widetilde{\psi}(T_\psi^n \dot{x}),
\]	
where $\dot{x}$ is the projection of $x\in X$ in the quotient group $X/C_\psi^\perp$ and $\tilde{\psi}$ is the character on $X/C_\psi^\perp$ induced by $\psi$. In particular, the observable $\psi(T^n x)$ on  $\mathscr{X} = (X, T)$ is equal to $\tilde{\psi}(T_{\psi}^n \dot{x})$ on $\mathcal{X}_\psi = (T_\psi, X/C_\psi^\perp)$, where as we have shown $X/C_\psi^\perp \cong \mathbb{T}^r \times C$ for some $r \geqslant 0$ and finite abelian group $C$.  Thus it suffices to prove \eqref{thm1goal} for all $x\in \mathbb{T}^r \times C$ and $\widetilde{\psi}\in \widehat{\mathbb{T}^r \times C}$ with $T$ be of the form \eqref{eqT}. 
This completes the reduction to the toral case.

\subsection{Affine linear maps on a torus.}

We have reduced Theorem 1.1 to the case that $\mathscr{X} = (T, X)$ with $X = \mathbb{T}^r \times C$ with $C$ finite and $T$ in the form \eqref{eqT}
 and {\color{black}distal}. For our purpose of examining observables $\xi(n)$ in this flow, we can ``linearize'' the flow by doubling the number of variables. That is, we will consider $Y = X \times X$ and the linear automorphism $W$ on $Y$ defined by
	\[
	W(x_1, x_2) = (A x_1 + x_2, x_2). 
	\]
Clearly, $\mathscr{Y} = (Y, W)$ is {\color{black}distal} since $\mathscr{X}$ is so, and the orbit $\{W^n(x_1, b)\}_{n \geqslant 1}$ is equal to $\{(T^n x_1, b)\}_{n \geqslant 1}$. Hence it suffices to prove the theorem for such $\mathscr{Y}$'s. That is we can assume that $\mathscr{X} = (X, W)$ with $X = \mathbb{T}^m \times F$, $F$ finite and $W$ a linear {\color{black}distal} automorphism of $X$. 
	
%Either by noting that the induced action of $W$ on $\widehat{X}$ must preserve $\{0\} \times \hat{F}$ (since these are precisely the elements of finite order in $\widehat{X}$) or 
The continuity of $W$ implies that it preserves the connected component $\mathbb{T}^m \times \{0\}$ of $0$ in $X$. Then $W$ must take the block triangular form
	\[
	W(\theta, f) = (B\theta + C f, D f),  
	\]
where $B: \mathbb{T}^m \to \mathbb{T}^m$ is an automorphism of the (connected) torus $\mathbb{T}^m$, $C : F \to \mathbb{T}^m$ is a homomorphism and $D : F \to F$ is an automorphism of $F$. The automorphism $B$ lifts to a linear automorphism $\widetilde{B}$ of $\mathbb{R}^m$ which preserves $\mathbb{Z}^m$, so that $\widetilde{B} \in GL_m(\mathbb{Z})$. {\color{black}Taking $f=0$ and noticing that  $W$ is distal, we see that $B$ is also distal}, and it is known that this implies that $\widetilde{B}$ is quasi-unipotent \cite{Dani}. That is,  $\widetilde{B}^{\nu_1} = U$ is unipotent for some $\nu_1 \geqslant 1$, i.e., 
$$
\widetilde{B}^{\nu_1} = U = I + N_1
$$
with $N_1$ nilpotent and $I$ the identity matrix. Also, since $F$ is finite, there is some $\nu_2 \geqslant 1$ such that $D^{\nu_2} = I$. Denote $ \text{lcm}(\nu_1, \nu_2)=\nu$. Then we have
	\[
	W^{\nu}(\theta, f) = \big((I + N_2)\theta + C_1 f, f\big), 
	\]
where  $N_2$ is nilpotent with $N_2^k=0$  for some $k \geqslant 0$ and $C_1$ is a morphism from $F$ to $\mathbb{T}^m$. Write
	\[
	\Phi := W^{\nu} = I + N,
	\]
where $N: X \to X$ is defined as $N(\theta,f)=(N_2\theta+C_1f, 0)$ satisfying $N^{k+1} \equiv 0$. Thus for a given integer $q \geqslant 0$ we have
	\[
	\Phi^q = \sum_{t=0}^{q} \binom{q}{t} N^t = \sum_{t=0}^{\min(k,q)} \binom{q}{t} N^t. 
	\]
	Writing $n \geqslant 0$ as $n = q\nu + l$ with $0 \leqslant  l < \nu$ we have
	\[
	W^n = W^{q\nu + l} = \Phi^q W^l.
	\] 
Hence, for any $x \in X$,  we have
\[
W^n x = \sum_{t=0}^{\min(k,q)} \binom{q}{t} N^t W^l x = \sum_{t=0}^{\min(k,q)} \binom{q}{t} \xi_{l,t}, 
\]
where
\[
\xi_{l,t} := \xi_{l,t} (x)= N^t W^l x. 
\]
For $\psi \in \widehat{X}$ fixed and $q$ varying with $q \geqslant k$, we have
\begin{equation} \label{eqWx}
	\psi(W^{q\nu+l} x) 
	= 
	\psi\left(\sum_{t=0}^{k} \binom{q}{t} \xi_{l,t}\right) 
	= \psi(\xi_{l,0})\psi(\xi_{l,1})^{\binom{q}{1}}\psi(\xi_{l,2})^{\binom{q}{2}} \cdots \psi(\xi_{l,k})^{\binom{q}{k}}.  
\end{equation}
The character $\psi \in \widehat{X}$ has the form $\psi : x \mapsto {\textup{e}}(\langle v, x \rangle)$ for some $v = (v_1, \dots, v_m) \in \mathbb{Z}^m$ where $\langle v, x \rangle$ means the dot product in $\mathbb{R}^m$. Write $\psi(\xi_{l,j})=\beta_{l,j}$ for $j=0,dots,k$. Then we can write \eqref{eqWx} as 
$$
\psi(W^{q\nu+l} x) 
= 
{\textup{e}}\left(\sum_{t=0}^{k} \binom{q}{t} \beta_{l,t}\right) 
$$
and hence the right-hand side of \eqref{eqWx} is ${\textup{e}}(Y(q))$, where $Y(q)$ is a polynomial in $q$ with degree $\leqslant  k$ and with coefficients depending on $v$ and the $\xi$'s. Changing variables from $q$ to $n$ by $n = \nu q + l$ with $0 \leqslant  l \leqslant  \nu-1$, we see that $Y(q) = P(n)$ is a polynomial in $n$ with degree $\leqslant  k$ and coefficients depending on $v, \nu, l$ and the $\xi$'s. It follows that
\begin{align*}
	\sum_{N<n\leqslant  N+M} \mu(n)\psi(W^n x) &= \sum_{l=0}^{\nu-1} \sum_{\substack{N<n\leqslant  N+M \\ n \equiv l \pmod{\nu}}} \mu(n)\psi(W^n x) \\
	&= \sum_{l=0}^{\nu-1} \sum_{\substack{N<n\leqslant  N+M \\ n \equiv l \pmod{\nu}}} \mu(n){\textup{e}}(P(n)).  
\end{align*}
Theorem \ref{thm1} now follows from the new  result of  Matom{\"a}ki et al \cite{TerenceTao2023Higheruniformity}  on  the estimation of exponential sum over primes in short intervals in Lemma \ref{youlideguji}.

\section{Proof of Theorem \ref{thm2} } 
	Since $A$ is a group automorphism of $\mathbb{T}^2$, we may take $A$ as a matrix in $GL_2(\mathbb{Z})$
By Lemma \ref{lemh}, the entropy of $\mathscr{X}$ is 
$0=h_{\text{top}}(T)=\sum_{|\lambda_i|>1}\log|\lambda_i|$, where $\lambda_i$ is any eigenvalue of $A$. Hence we have $|\lambda_i|\leqslant 1$ for each eigenvalue of $A$.  

Write $t_A=\operatorname{tr}(A), d_A=\det(A)=\pm1$ and $\Delta = t_A^2-4d_A$. Then the characteristic equation and the discriminant of $A$ are
$$
\lambda^2-t_A\lambda+d_A=0,
\qquad
\Delta = t_A^2-4d_A
$$
respectively.
\begin{enumerate}
 \item  If $\Delta>0$, then there are two distinct real eigenvalues. Notice that their product is $\pm 1$. If one of them has modulus larger than $1$, then the topological entropy is positive, a contradiction. Hence both of the two distinct real eigenvalues must have absolute value $1$, i.e. $
 \lambda_1=1$, $\lambda_2=-1$, which can be viewed as unit roots of order $2$.  
 
 \item  If $\Delta<0$, then $d_A=1$, otherwise  $d=-1$ will imply $\Delta = t^2+4>0$, a contradiction. Thus $d_A=1$, which combined with $\Delta<0$ implies $t_A^2<4$, so $t_A=0,\pm1$. From the  characteristic equation  $\lambda^2-t_A\lambda+d_A=0$ we get $\lambda=\pm i$ or $\lambda=\frac{1}{2}(1\pm\sqrt{3}i)$,
 or 
 $\lambda=\frac{1}{2}(-1\pm\sqrt{3}i)$,  which implies that the eigenvalues of $A$ are unit roots of order $3,4,6$. In this case $\mathscr{X}$ is of  zero entropy
 
 \item  If $\Delta=0$, then $t_A^2=4d_A$. Since $d_A=\pm1$, we get $d_A=1$ and $t_A=\pm 2$. Then the eigenvalue is a repeated root $1$ or $-1$, and the entropy is zero. In this case the eigenvalue can also be viewed as unit roots of order $2$. 
\end{enumerate} 
In summarize,  the eigenvalues of $A$ are unit roots of order $2,3,4,6$. 

For any natural numbers $0\leqslant  k\leqslant  n$, denote
$$
A^k=
\begin{pmatrix}
	a_1(k) & a_2(k) \\
	a_3(k) & a_4(k)
\end{pmatrix},
\qquad
\sum_{k=0}^{n-1}A^k b=
\begin{pmatrix}
	d_1(n)   \\
	d_2(n) 
\end{pmatrix}.
$$
Then applying \eqref{eqT} repeatedly we get that for any $x=(x_1,x_2)\in \mathbb{T}^2$,
$$
T^n(x)
=
A^nx+\sum_{k=0}^{n-1}A^kb
=
\Big(a_1(n)x_1 + a_2(n)x_2 +d_1(n),
a_3(n)x_1 + a_4(n)x_2 +d_2(n)\Big)	 
$$
Then, for any $c=(c_1,c_2)\in \mathbb{Z}^2$, we have 
\begin{align*}
\langle c,T^n(x)\rangle
&=
c_1\big(a_1(n)x_1 + a_2(n)x_2 +d_1(n)\big)
+
c_2\big(a_3(n)x_1 + a_4(n)x_2 +d_2(n)\big)
\\
%&=c_1 a_1(n)x_1 + c_1a_2(n)x_2 +c_1d_1(n) 
%+
%c_2a_3(n)x_1 + c_2a_4(n)x_2 + c_2d_2(n)
%\\
&=\big(c_1 a_1(n) + c_2a_3(n)\big)x_1  
+ \big(c_1a_2(n) + c_2a_4(n)\big) x_2 
+c_1d_1(n) +c_2d_2(n).
\end{align*}
It follows that for any $x\in X$ 
\begin{align*}\label{smn}
&	S(M,N)
	:=\sum_{N< n\leqslant  N+M}\mu(n) {\textup{e}}(\langle c, T^n(x)\rangle) \notag
	\\
	&\ll
	\sum_{N<n\leqslant  N+M}\mu(n) 
	{\textup{e}}\Big(
	\big(c_1 a_1(n) + c_2a_3(n)\big)x_1  
	+ \big(c_1a_2(n) + c_2a_4(n)\big) x_2 
	+c_1d_1(n) +c_2d_2(n)
	\Big).
\end{align*}

By virtue of the Jordan canonical form of 2‑by‑2 matrices, $A^n$ takes one of the following two forms: 
$$
Q
\begin{pmatrix}
	\lambda^n & n \lambda^{n-1}\\
	0 & \lambda^n
\end{pmatrix} 
Q^{-1},
\qquad
Q
\begin{pmatrix}
	\lambda^n_1& 0 \\
	0& \lambda_2^n
\end{pmatrix}
Q^{-1},
$$
where $Q$ is an invertible matrix and $\lambda$, $\lambda_i$ are eigenvalues of $A$, which are unit roots of order $2,3,4,6$. 
Then for each fixed residue class $n \equiv k$ (mod 12), the entries $a_i(n)$ of $A^n$ are polynomials in $n$ of degree at most $1$, and  the entries $d_i(n)$ of the translation terms are polynomials in $n$ of degree at most $2$.
Therefore, applying Lemma \ref{youlideguji} we get
$$
S(M,N)
\ll
\sum_{k=1}^{12}\sum_
{\substack
	{N < n \leqslant  N+M
\\
n\equiv k\pmod{12}}
}
 \mu(n)\, {\textup{e}}\big(P(n) \big)
\ll
\frac{M}{\log^{B}N}
$$ 
for any  $B>0$, where $P(n)$ is a polynomial in $n$ of degree no more than $2$.
 
Finally, applying Jackson's theorem, any Lipschitz functions on $\mathbb{T}^2$ can be approximated by linear combinations of finitely many trigonometric polynomials. Therefore, the same upper bound also holds for Lipschitz functions on $\mathbb{T}^2$. This complete the proof.

\bigskip
\noindent{\bf{Acknowledgements}} Liu is supported by the National Natural Science 
Foundation of China (Grant No. 12401012).

\end{document}